%% file: msri2-.tex
\documentstyle[11pt,mssymb]{article}
\def\mra{\rightarrow}
\def\ra{$\mra$}
\def\mlra{\longrightarrow}
\def\lra{$\mlra$}
\def\mex{\mathop{\rm mex}}
\def\nimsum{\stackrel{*}{\textstyle +}}
\def\ul{\underline}
\def\0{^{\phantom0}}
\def\sb{\_} % stopped blocked file
\input{macros2}
\begin{document}

\vspace*{-6ex}

\centerline
{\Large\bf Higher Nimbers in pawn endgames on large chessboards}

\vspace*{4ex}

\centerline{\large Noam D. Elkies}

\vspace*{5ex}

{\small
{\em
Do $*2$, $*4$ and higher Nimbers occur
on the $8\times 8$ or larger boards?
}
\\
--- ONAE [Elkies 1996, p.148]

\vspace*{2ex}

{\em It's full of stars!} 
\\
--- 2001: A Space Odyssey [Clarke 1968, p.193]
}

\begin{quote}
{\small
{\sc Abstract.} We answer a question posed in [Elkies~1996]
by constructing a class of pawn endgames on $m\times n$ boards
that show the Nimbers $*k$\/ for many large~$k$.  We do this by
modifying and generalizing T.R.~Dawson's ``pawns game''
[Berlekamp et al.~1982].  Our construction works for $m\geq 9$
and $n$ sufficiently large; on the basis of computational evidence
we conjecture, but cannot yet prove, that the construction yields
$*k$ for all integers~$k$.
}
\end{quote}

\centerline{ {\large\bf 1.~Introduction} }

In [Elkies 1996] we showed that certain chess endgames can be
analyzed by the techniques of combinatorial game theory (CGT).
We exhibited such endgames whose components show a variety of
CGT values, including integers, fractions, and some infinite
and infinitesimal values.  Conspicuously absent were the values
$*k$\/ of Nim-heaps of size $k>1$.  Towards the end of~[Elkies 1996]
we asked whether any chess endgames, whether on the standard
$8\times 8$ chessboard or on larger rectangular boards, have
components equivalent to $*2$, $*4$ and higher Nimbers.
In the present paper we answer this question affirmatively
by constructing a new class of pawn endgames on large boards
that include $*k$ for many large~$k$, and conjecture --- though we
cannot yet prove --- that all $*k$ arise in this class.

Our construction begins with a variation of a pawns game called
``Dawson's Chess'' in~[Berlekamp et al.~1982].  In~\S3 we introduce
this variation and show that, perhaps surprisingly, all quiescent
(non-entailing) components of the modified game are equivalent to
Nim-heaps (Thm.~1).  We then determine the value of each such
component, and characterize non-loony moves (Thm.~2).  In~\S4
we construct pawn endgames\footnote{
  More precisely, King-and-pawn endgames; the first two words are
  usually suppressed because every legal chess position must have
  a White and a Black king.
  }
on large chessboards that incorporate those components.
These endgames do not yet attain our aim, because
the values determined in Thm.~2 are all $0$ or~$*1$.
In~\S5 we modify one of our components to obtain $*2$.  
In~\S6 we further study components modified in this way,
showing that they, too, are equivalent to Nim-heaps (Thm.~3).
We conclude with the numerical evidence suggesting
that all Nim-heaps can be simulated by components
of pawn endgames on large rectangular chessboards.

Before embarking on this course, we show in~\S2 a pair of
endgames on the standard $8\times 8$ board in the style of
[Elkies~1996] that illustrate the main ideas, specifically
the role of ``loony moves''$\!$.  Readers who are
much more conversant with CGT than with chess endgames
will likely prefer to skim or skip \S2 on first reading,
returning to it only after absorbing the theory in~\S3.
Conversely, chessplayers not fluent in~CGT will find in~\S2
motivation for the CGT ideas central to \S3 and later sections.

\vspace*{3ex}

\centerline{ {\large\bf 2.~An illustrative pair of endgames} }

We introduce the main ideas of our construction by analyzing
the following pair of composed endgames on
the standard $8\times 8$ chessboard:\footnote{
  These are not pawn endgames, but all units
  other than pawns are involved in the Kingside Zugzwang,
  and are thus passive throughout the analysis.  We can
  construct plausible positions where that Zugzwang, too,
  is replaced by one using only Kings and pawns, but only at the cost
  of introducing an inordinate number of side-variations tangential
  to the CGT content of the positions.  For instance, replace files
  e--g in Diag.~1A by White Kg1, Ph2 and Black Kh3, Pg2; in Diag.~1B,
  do the same and move the h5/h7 pawns to e3/e5.
  }

\vspace*{3ex}

\centerline{
\board
{........}
{........}
{ppp.....}
{........}
{PPP.....}
{....k.b.}
{....pNP.}
{....K.Q.}
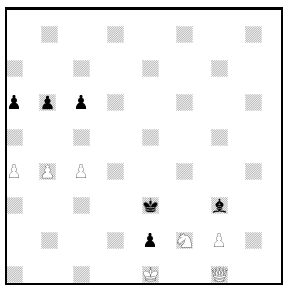
{Diag.~1A: whoever moves loses}
\qquad\quad
\board
{........}
{.......p}
{ppp.....}
{.......P}
{PPP.....}
{....k.b.}
{....pNP.}
{....K.Q.}
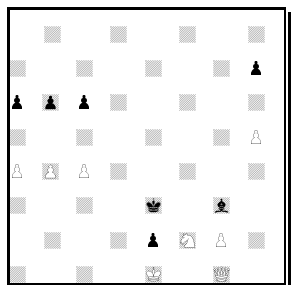
{Diag.~1B: whoever moves wins}
\
}

Diagram~1A consists of two components.  On the Kingside, seven men
are locked in a mutual Zugzwang that we already used in [Elkies~1996].
Both sides can legally move in the Kingside, but only at the cost of
checkmate (Qh1(2) Bxf2) or ruinous material inferiority.
Thus the outcome of Diagram~1A hinges on the Queenside component,
with three adjacent pawns on each side.
We have not seen such a component in [Elkies~1996],
but it turns out to be a mutual Zugzwang: whichever side moves first,
the opponent can maneuver to make the last pawn move on the Queenside,
forcing a losing King move in the Kingside component.  Thus a5 can
be answered by bxa5 bxa5 c5, likewise c5 by bxc5 bxc5 a5, and b5
by axb5 axb5 c5.  In this last line, it is no better to answer
axb5 with cxb5, since then c5 wins: if played by Black, White will
respond a5 and promote first, but Black c1Q will be checkmate;
while if White plays c5 and Black answers a5, White replies c6
and promotes with Black's pawn still two moves away from the first rank,
winning easily.\footnote{
  This part of the analysis explains why we chose this Kingside
  position from [Elkies~1996]: the position of White's King,
  but not Black's, on its first rank makes it more vulnerable
  to promoted pawns, exactly compensating for the White pawns
  being a step closer to promotion than Black's.  In [Elkies~1996]
  this Kingside Zugzwang had the Kings on f1 and f3, not e1 and e3;
  here we shifted the position one square to the left
  so as not to worry about a possible White check by a newly promoted
  Queen on a8.  Thus 1~b5 can be answered by 1\ldots cxb5
  as well as 1\ldots axb5.
  }
Note the key point that a5 or c5 must not be answered by b5?, since
then c5 (resp.~a5) transfers the turn to the second player and wins
--- but not cxb5?\ cxb5 (or axb5?\ axb5), regaining the Zugzwang.

Diagram~1B is Diag.~1A with the h5 and h7 pawns added;
these form an extra component which we recognize from [Elkies~1996]
as having the value~$*$.  Thus we expect that Diag.~1B is a
first-player win, and indeed either side can win by playing h6,
effectively reducing to Diag.~1A.  The first player can also win
starting on the Queenside:
a5 bxa5 bxa5 ($*+*=0$), and likewise if a5 is answered by b5
(c:b5 etc.\ wins, but not c5?\ h6).  The first player must
not, however, start b5?, when the opponent trades twice on b5,
in effect transforming $*+*$ into $0+*$, and then wins by
playing h6.

Note the role of the move b5 in the analysis of both Diagrams 1A
and~1B.  In the terminology of [Berlekamp et al.~1982], this is
an ``entailing move'': it makes a threat (to win by capturing
a pawn) that must be answered in the same component.  But, whether
the rest of the position has value~$0$ (Diag.~1A) or * (Diag.~1B),
the move b5 loses, because the opponent can answer the threat
in two ways, one of which passes the turn back (advancing
the threatened pawn to a5 or c5), one of which in effect
retains the turn (capturing on b5, then making another move
after the forced re-capture).  Since the first option wins if the
rest of the position has value~$0$, and the second wins if the
rest of the position has value~$*$ (or any other nonzero Nimber),
b5 is a bad move in either case.  In [Berlekamp et al.~1982]
such an unconditionally bad move is called ``loony'' (see p.378).
Since b5 is bad, it follows in turn that, in either Diag.~1A or~1B,
the entailing move a5 may be regarded as the non-entailing ``move''
consisting of the sequence a5 bxa5 bxa5: the only other reply to~a5
is~b5, which is loony and so can be ignored.

We next show that this analysis can be extended to similar
pawn configurations on more than three adjacent files.
We begin our investigation with a simplified game
involving only the relevant pawns.

\vspace*{3ex}

\centerline{ {\large\bf 3.~A game of pawns} }

{\bf 3.1: Game definition}

Our game is played on a board of arbitrary length~$n$ and height~$3$.
At the start, White pawns occupy some of the squares on the bottom row,
and Black pawns occupy the corresponding squares on the top row.
For instance, Diag.~2 shows a possible starting position
on a board with $n=12$:

\vspace*{2ex}

\centerline{
\longboardthree
{ppppppp.pppp}
{............}
{PPPPPPP.PPPP}
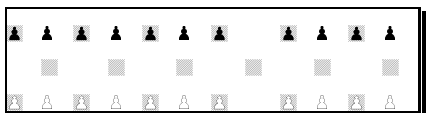
{Diag.~2: A starting position}
\
}

The pawns move and capture like chess pawns, except that there is
no double-move option (and thus no {\em en passant}\/ rule).  Thus
if a file contains a White pawn in the bottom row and a Black pawn
in the top row, these pawns were there in the initial position and
have not moved; we call such a file an ``initial file''$\!$.

White wins if a White pawn reaches the top row, and Black wins if
a Black pawn reaches the bottom row.  Thus there is no need for a
promotion rule because whoever could promote a pawn immediately wins
the game.  But it is easy enough to prevent this, and we shall assume
that neither side allows an opposing pawn to reach its winning row.
The game will then end in finitely many moves with all pawns blocked,
at which point the winner is the player who made the last move.
We shall sometimes call this outcome a ``win by Zugzwang''$\!$,
as opposed to an ``immediate win'' by reaching the opposite row.

For instance, from Diag.~2 White may begin 1~c2.
Since this threatens to win next move by capturing on b3 or d3,
Black must capture the c2-pawn; if Black captures with the b-pawn,
we reach Diag.~3A.  Now Black threatens to win by advancing this
pawn further, so White must capture it; but unlike Black's capture,
White's can only be made in one way: if 2~dxc2?, threatening to win
with 3~cxd3, Black does not re-capture but instead plays d2, producing
Diag.~3B.  Black then wins, since touchdown at d1 can only be delayed
by one move with 3~exd2 exd2 (but not 3\ldots cxd2??, when 4~c3 wins
instead for White!), and then 4\ldots d1.

\vspace*{2ex}

\centerline
{
\longboardthree
{p.ppppp.pppp}
{..p.........}
{PP.PPPP.PPPP}
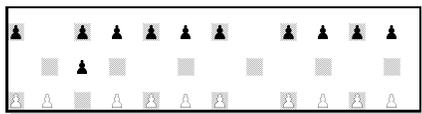
{Diag.~3A: After 1 c2 bxc2}
\qquad
\longboardthree
{p.p.ppp.pppp}
{..Pp........}
{PP..PPP.PPPP}
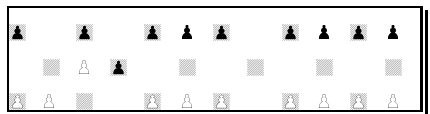
{Diag.~3B: If then 2 dxc2?\ d2, winning}
\
}

Thus White must play 2~bxc2 from Diag.~3A.  This again threatens to
win with 3~cxd3, so Black has only two options.  One is to re-capture
with 2\ldots dxc2, forcing White in turn to re-capture: 3~dxc2,
reaching Diag.~3C.  Alternatively, Black may save the d3-pawn by
advancing it to~d2 (Diag.~3D).  This forces White to move the attacked
pawn at e1.  Again there are two options.  White may capture with
3~exd2, forcing Black to reply 3\ldots exd2 (Diag.~3E), not cxd2?\
when White wins immediately with 4~c3.  Alternatively White may
advance with 3~e2 (Diag.~3F), when Black again has two options
against 4~exf3, etc.  Eventually the skirmish ends either with
mutual pawn captures (as in Diag.~3C or~3E) or when the wave
of pawn advances reaches the end of the component (3\ldots f2 4~g2),
leaving one side or the other to choose the next component to play in.

\vspace*{2ex}

\centerline
{
\longboardthree
{p.p.ppp.pppp}
{..P.........}
{P...PPP.PPPP}
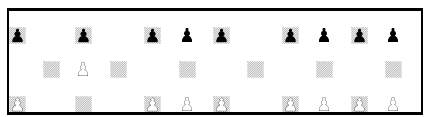
{Diag.~3C: After 2 bxc2 bxc2 3 dxc2}
\quad
\longboardthree
{p.p.ppp.pppp}
{..Pp........}
{P..PPPP.PPPP}
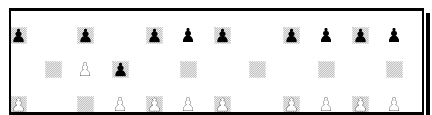
{Diag.~3D: Instead 2\ldots d2}
\
}

\vspace*{2ex}

\centerline
{
\longboardthree
{p.p..pp.pppp}
{..Pp........}
{P..P.PP.PPPP}
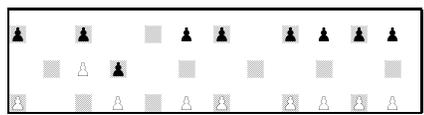
{Diag.~3E: Then 3~exd2 exd2}
\quad
\longboardthree
{p.p.ppp.pppp}
{..PpP.......}
{P..P.PP.PPPP}
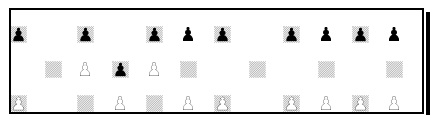
{Diag.~3F: or 3~e2}
\
}

As noted in the Introduction, this game is reminiscent of
the game called ``Dawson's Chess'' in [Berlekamp et al.~1982,
pp.~88~ff.]; the only difference is that in Dawson's Chess a player
who can capture a pawn must do so, while in our game, as in ordinary
chess, captures are optional.\footnote{
  Since Dawson was a chess problemist, we first guessed that
  the game analyzed here was Dawson's original proposal, before the
  modification in [Berlekamp et al.~1982].  But in fact R.K.~Guy
  reports in a 16.viii.1996 e-mail that Dawson did want obligatory
  captures but proposed a mis\`ere rule (last player loses).  Guy
  also writes: ``I am aware of some very desultory attempts to analyze
  the game in which captures are allowed, but little was achieved,
  to my knowledge.''  I thank Guy for this information, and John
  Beasley who more recently sent me copies of pages from the 12/1934
  and 2/1935 issues of {\em The Problemist}\/ Fairy Supplement in which
  Dawson proposed and analyzed his original game.
  }
Because of the obligation to capture,
Dawson's Chess may appear to be an entailing game, but it is quickly
seen to be equivalent to a (non-entailing) impartial game, called
``Dawson's Kayles'' in [Berlekamp et al.~1982].  Our game also has
entailing moves, and features a greater variety of possible components;
but we shall see that it, too, reduces to a non-entailing impartial
game once we eliminate moves that lose immediately and loony moves.

{\bf 3.2: Decomposition into components}

We begin by listing the possible components.  We may ignore any
files in which no further move may be made.  These are the
files that are either empty (such as the h-file throughout Diag.~3,
the b-file in Diag.~3C--3F, and the d- and e-files in 3C and~3E
respectively) or closed.  We say a file is ``closed'' if it contains
one pawn of each color, neither of which can move or capture
(the c-file in 3C--3F and the d-file in 3E--3F).  One might object
that such a currently immobile pawn may be activated in the future;
for instance in Diagram~3D if White plays 3~exd2 then the dormant c-file
may awake: 3\ldots cxd2.  But we observed already that this Black move
loses immediately to 4~c3.  Since we may and do assume that
immediately losing moves are never played, we may ignore the
possibility of 3\ldots cxd2?, and regard the c-file in Diag.~3D
and the d-file in Diag.~3F as permanently closed.

By discarding empty or closed files we partition the board into
components that do not interact except when an entailing move
requires an answer in the same component.  Thus at each point
there can be at most one entailing component (again assuming no
immediately losing moves).  We next describe all possible components
and introduce a notation for each.

A non-entailing component consists of $m$ consecutive initial files,
for some positive integer~$m$.  We denote such a component by~$[m]$.
For instance, Diag.~3A is $[7]+[4]$; Diag.~3C is $[1]+[3]+[4]$;
and Diag.~3E is $[1]+[2]+[4]$.  An entailing component contains
a pawn that has just moved (either vertically or diagonally)
to the second rank, threatening an immediate win, and can be captured.
We denote this pawn's file by ``:''$\!$.
There are four kinds of entailing component, depending on whether
this pawn can be captured in one or two ways
and on how many friendly pawns defend it
by being in position to re-capture.

$\bullet$ If the pawn is attacked once and defended once,
then the attacking and defending pawns are on the same file
(else the side to move can win immediately as in Diag.~3B).
Thus the ``:''~file is at the end of a component
each of whose remaining files is initial.
We denote the component by~[:$m$\/], where $m$ is the number
of initial files in the component.  For instance, Diag.~4A shows [:5].
For the mirror-image of~``[:$m$\/]''$\!$,
we use either the same notation or ``[$m$\/:]''.
Either side can move from $[m+1]$ to~[:$m$\/]
by moving the left- or rightmost pawn.
Faced with [:$m$\/], one must move the attacked pawn, either by
capturing the ``:''~pawn or by advancing it.  Advancing yields
[:$(m-1)$] (see Diag.~4B), unless $m=1$ when the advance yields~$0$
since all files in the component become blocked.
Capturing yields [:.]$+[m-1]$ (see Diag.~4C),
where [:.] is the the component defined next.

$\bullet$ If the pawn is attacked once and not defended, then
it has just captured, and is subject to capture from an unopposed pawn.
We denote the file with one unopposed pawn by~``.''.
The capture is obligatory, and results in a closed file.
Therefore the adjacent ``:''\ and~``.'' files to not interact
with any other components, even if they are not yet separated from them
by empty or closed files.  We may thus regard these files as
a separate component~[:.], which entails a move to~$0$.  For
instance, files $b,c$ in Diag.~4C constitute [:.],
which is unaffected by the $[4]$ on files d through~g.

\vspace*{2ex}

\centerline
{
\boardthree
{..ppppp.}
{.p......}
{.PPPPPP.}
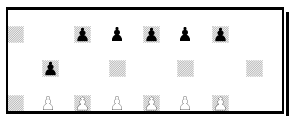
{Diag.~4A: [:5]}
\
\boardthree
{..ppppp.}
{.pP.....}
{.P.PPPP.}
{Diag.~4B: [:4]}
\
\boardthree
{..ppppp.}
{.P......}
{.P.PPPP.}
{Diag.~4C: [:.]+[4]}
\
}

$\bullet$ If the pawn is attacked twice and defended twice,
then the component of the ``:''~file consists of that file,
$m$ initial files to its left, and $m'$ initial files to its right,
for some positive integers $m$ and~$m'$.
We denote such a component by~[$m$\/:$m'$].
For instance, Diag.~5A shows [2:4].
We already encountered this component after the move c2 from Diag.~2.
The component [$m$\/:$m'$] entails a capture of the ``:''~pawn.  
This yields either $[m-1]+$[.:$m'$] or [.:$m]+[m'-1]$, where
[.:$m\/$], defined next, is our fourth and last kind of entailing
component, and $[m-1]$ (or $[m'-1]$) is read as~$0$ if $m=1$
(resp.~$m'=1$).

$\bullet$ If the pawn is attacked twice and defended once,
then it is part of a component obtained from [:.] by placing
$m$ initial files next to the ``:''~file, for some positive
integer~$m$.  Naturally we call such a component [.:$m$\/]
(or [$m$\/:.]).  As explained in the paragraph introducing [:.],
an initial file next to the ``.''~file cannot interact with it,
and thus belongs to a different component.  For instance, the
two possible captures from [2:4] yield [1]+[.:4] (Diag.~5B,
also seen in Diag.~3A) and [.:2]+[3] (Diag.~5C).  Faced with
[.:$m$\/], one has a single move that does not lose immediately:
capture with the ``.''~pawn, producing [:$m$\/], as seen earlier
in connection with Diag.~3C.

\vspace*{2ex}

\centerline
{
\boardthree
{ppppppp.}
{..P.....}
{PP.PPPP.}
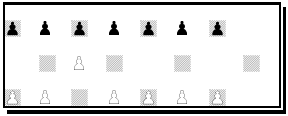
{Diag.~5A: [2:4]}
\
\boardthree
{p.ppppp.}
{..p.....}
{PP.PPPP.}
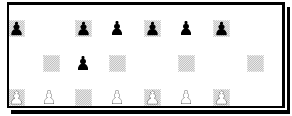
{Diag.~5B: [1]+[.:4]}
\
\boardthree
{ppp.ppp.}
{..p.....}
{PP.PPPP.}
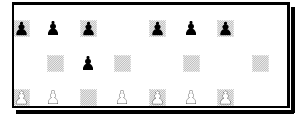
{Diag.~5C: [.:2]+[3]}
\
}

We summarize the available moves as follows.  It will be convenient
to extend the notations $[m]$, [:$m$\/], [.:$m$\/], [$m$\/:$m'$]
by allowing $m=0$ or $m'=0$, with the understanding that
$$
[0]=[:\!0]=[0\!:\!0]=0, \quad [.\!:\!0]=[:\!.], \quad
[m\!:\!0]=[0\!:\!m]=[:\!m].
$$
We then have:
\begin{itemize}
\item From $[m]$, either side may move to
$[m_1\!:\!m_2]$ for each $m_1,m_2\geq0$ such that $m_1+m_2=m-1$.
\item If $m>0$ then $[:\!m]$ entails a move to
either $[:\!.]+[m-1]$ or $[:\!(m-1)]$.
\item $[:\!.]$ entails a move to~$0$.
\item If $m,m'>0$ then $[m\!:\!m']$ entails a move to either
$[m-1]+[.\!:\!m']$ or $[m'-1]+[.\!:\!m]$.
\item If $m>0$ then $[.\!:\!m]$ entails a move to $[:\!m]$.
\end{itemize}
For instance, Diag.~3D shows $[1]+[:\!3]+[4]$, moving either
to $[1]+[:\!.]+[2]+[4]$ and thence to $[1]+[2]+[4]$ (Diag.~3E),
or to $[1]+[:\!2]+[4]$ (Diag.~3F).

Our list of possible moves confirms that
$[m]$, [:$m$\/], [:.], [.:$m$\/], [$m$\/:$m'$]
are the only possible components: the initial position is
a sum of components $[m_i]$, and each move from a known component
that does not concede an immediate win
yields a sum of $0$, $1$, or $2$ known components.

{\bf 3.3: Analysis of components}

Since in each component both sides have the same options,
we are dealing with an {\em impartial}\/ entailing game.
We could thus apply the theory of such games, explained
in~[Berlekamp et al.~1982], to analyze each component.
But it turns out that once we eliminate loony moves
the game is equivalent to an ordinary impartial game,
and thus that each component $[m]$ is equivalent to
a Nim-heap of size depending on~$m$.  

Consider the first few~$m$.
Clearly $[1]$ is equivalent to $*1$, a Nim-heap of size~$1$.
At the end of \S2 we have already seen in effect
that a move from~[2] to~[:1] is loony: if the rest of the
game has value~0, then the reply [:1]\ra0 wins;
otherwise, the reply [:1]\ra[:.] forces [:.]\ra0,
again leaving a forced win in the sum of the remaining components.
Thus also [3]\ra[1:1] is loony, because the forced continuation
[1:1]\ra[.:1]\ra[:1] again leaves the opponent in control.
On the other hand, [3]\ra[:2] is now seen to force
[:2]\ra[:.]+[1], since the alternative [:2]\ra[:1] is known to lose.
We thus have the forced combination [3]\ra[:2]\ra[:.]+[1]\ra[1],
which amounts to a ``move'' [3]\lra[1].\footnote{
  Here and later, we use an arrow ``\ra'' for a single move,
  and a long arrow ``\lra'' for a sequence of 3, 5, 7, \ldots\
  single entailing moves considered as one ``move''$\!$.
  }
Moreover, we have shown that this
is the only reasonable continuation from~[3].
Since \hbox{$[1]\cong *1$}, we conclude that
[3] is equivalent to a Nim-heap of size $\mex(\{1\})=0$,
i.e.\ a mutual Zugzwang, as we already discovered
in the analysis of Diags.~1A,1B.

What of~[4] and~[5]?  From~[4], there are again two options,
one of which can be eliminated because it leads to the loony~[:1],
namely [4]\ra[2:1] (after [2:1]\ra[1]+[.:1]\ra[1]+[:1]).
This leaves [4]\ra[:3], which in turn allows two responses.
One produces the sequence [4]\ra[:3]\ra[2]+[:.]\ra[2],
resulting in a value of~$0$.
The other response is [4]\ra[:3]\ra[:2], which we already know forces
the further [:2]\ra[:.]+[1]\ra[1]$\,\cong *1$.  In effect, the
response to [4]\ra[:3] can interpret the move either as
$[4]\mlra[2]\cong0$, or as $[4]\mlra[1]\cong *1$
{\em with the side who played $[4]\mra[:\!3]$ on move again}.
We claim that the latter option can be ignored.  The reason is that
the first interpretation wins unless the remaining components
of the game add to~$0$; but then the second interpretation
leaves the opponent to move in a nonzero position, and thus
also loses.  We conclude that [4] is equivalent to an impartial
game in which either side may move to~$0$, and is thus equivalent
to a Nim-heap of size $\mex(\{0\})=1$.  As to~[5], there are now
three options, only one excluded by~[:1], namely [5]\ra[3:1].
The option [5]\ra[2:2] (move the center pawn) forces the continuation
[2:2]\ra[1]+[.:2]\ra[1]+[:2]\ra[1]+[:.]+[1]\ra[1]+[1], and is thus
tantamount to $[5]\mlra[1]+[1]\cong0$.  This leaves [5]\ra[:4],
which we show is loony for a new reason.  One continuation is
[:4]\ra[:.]+[3]\ra[3], interpreting the move as $[5]\mlra[3]\cong0$.
The other is [:4]\ra[:3], which as we have just seen is equivalent to
[:4]\lra[2].  Since $[2]\cong0$, this continuation interprets [5]\ra[:4]
as a move to~$0$ {\em followed by an extra move}.  Thus a move to~[:4]
always allows a winning reply, namely [:4]\ra[:3] if the remaining
components add to~$0$, and [:4]\ra[:.]+[3] if not.  Hence the move
to~[:4] is loony as claimed, and $[5]\cong *(\mex\{0\}) = *1$.

The analysis of [1] through~[5] shows almost all the possible behaviors
in our game; continuing by induction we prove:

{\bf Theorem 1.}
{\sl
i) For each integer~$m\geq 0$, the component $[m]$
is equivalent to a Nim-heap of some size~$\epsilon(m)$.\\
ii) A move to $[:\!1]$ is loony.  For each integer $m>1$,
a move to $[:\!m]$ is either loony or equivalent to a move to
$[m-1]\cong *(\epsilon(m-1))$.  The move is loony if and only if
a move to $[:\!(m-1)]$ is \ul{not} loony and
$\epsilon(m-1)=\epsilon(m-2)$.\\
iii) For any positive integers $m_1,m_2$, a move to $[m_1\!:\!m_2]$
is either loony or equivalent to a move to
$[m_1-1]+[m_2-1] = *(\epsilon(m_1-1) \nimsum \epsilon(m_2-1))$.
The move to $[m_1\!:\!m_2]$ is non-loony if and only both~$m_1$
and $m_2$ satisfy the criterion of~(ii) for a move to~$[:\!m]$
to be non-loony.\\
iv) We have $\epsilon(0)=0$, $\epsilon(1)=1$,
and for $m>1$ the values $\epsilon(m)$ are given recursively by
$$
\epsilon(m) =
\mex_{m_1,m_2} \Bigl(\epsilon(m_1-1) \nimsum \epsilon(m_2-1)\Bigr).
$$
Here the mex runs over pairs $(m_1,m_2)$ of nonnegative integers 
such that $m_1+m_2=m-1$ and a move to $[m_1\!:\!m_2]$ is not loony,
as per the criteria in parts (ii) [for $m_1 m_2 = 0$]
and~(iii) [for $m_1 m_2 > 0$].
For this equation we declare that $\epsilon(-1)=0$.
}

In parts (iii) and (iv), ``$\nimsum$'' denotes the Nim sum:
\hbox{$(*k) + (*k') = *(k\nimsum k')$}.  Thus, once parts (i)--(iii)
are known for all components with fewer than~$m$ initial files,
(iv) is just the Sprague-Grundy recursion for impartial games.
Once (iv) is known for all $m\leq m_0$, so is~(i).
The arguments for~(ii) and~(iii) are the same ones we used
for components with up~to $5$ initial files.  For instance,
for (iii) the move $[m_1+1+m_2]\mra[m_1\!:\!m_2]$ forces a choice
among the continuations
[$m_1$:$m_2$]\ra[$m_1-1$]+[.:$m_2$]\ra[$m_1-1$]+[:$m_2$]
and 
[$m_1$:$m_2$]\ra[$m_1$:.]+[$m_2-1$]\ra[$m_1$:]+[$m_2-1$].
If a move to [:$m_1$] or [:$m_2$] is loony then one or both
of these continuations wins.  Otherwise by (ii) both continuations
are tantamount to $[m_1+1+m_2]\mlra[m_1-1]+[m_2-1]$.~~@p

Carrying out the recursion for $\epsilon(m)$,
we quickly detect and prove a periodicity:

{\bf Theorem 2.}
{\sl
For all $m\geq 0$ we have $\epsilon(m)=0$ if $m$ is congruent
to~$0$, $2$, $3$, $6$, or\/~$9\bmod 10$, and $\epsilon(m)=1$
otherwise.  A move to $[:\!m]$ is loony if and only if
$m=5k\pm1$ for some integer~$k$.
}

{\sl Proof}\/: Direct computation verifies the claim through $m=23$,
which suffices to prove it for all~$m$ as in [Berlekamp et al.~1982,
pp.~89--90], since $23$ is twice the period plus the maximum number
of initial files lost by a ``move'' $[m]\mlra[m-2]$ or
$[m_1+1+m_2]\mlra[m_1-1]+[m_2-1]$.~~@p

So, for instance, Diagram~2 is equivalent to
$*(\epsilon(7)) + *(\epsilon(4)) = *1 + *1 = 0$
and is thus a mutual Zugzwang, a.k.a.\ ${\cal P}$-position
or previous-player win.  The next player thus might as well
play a loony move such as 1~c2, in the hope of giving the
opponent Enough Rope [Berlekamp et al.~1982, p.17];
the only correct response is 1\ldots bxc2 (Diag.~3A)
2~bxc2 d2!\ (Diag.~3D), maintaining the win after either
3~exd2 exd2 (Diag.~3E) or 3~e2 (Diag~3F) fxe2!\ 4~fxe2.

\vspace*{3ex}

\centerline{ {\large\bf 4.~Embedding into generalized chess} }

Consider Diag.~6A, a pawn endgame on a chessboard of
$9$ rows by $12$ files:

\vspace*{2ex}

\centerline{
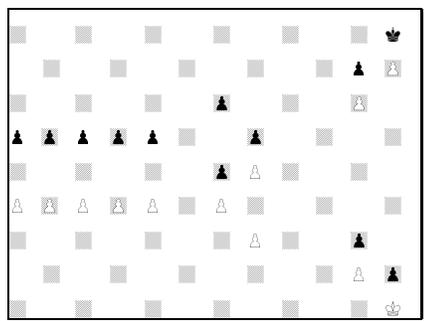
\def\caption{Diag.~6A: whoever moves wins (c5!)}
\longboardplus
{...........k}
{..........pP}
{......p...P.}
{ppppp..p....}
{......pP....}
{PPPPP.P.....}
{.......P..p.}
{..........Pp}
{...........K}
\quad
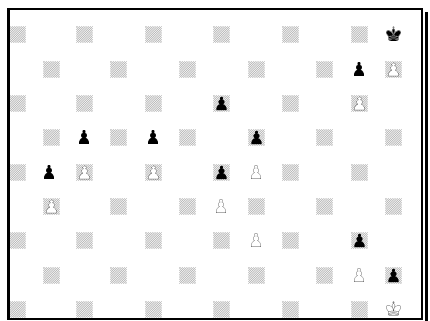
\def\caption{Diag.~6B: after 1 c5 \ldots\ 4 e5}
\longboardplus
{...........k}
{..........pP}
{......p...P.}
{..p.p..p....}
{.pP.P.pP....}
{.P....P.....}
{.......P..p.}
{..........Pp}
{...........K}
\
}

There are four components.
In each of the top right and bottom right corners,
a King and three pawns are immobilized.
Near the middle of the board (on the g- and h-files),
we have a mutual Zugzwang with three pawns on a side;
a player forced to move there will allow an opposing pawn
to capture and soon advance to Queen promotion, giving checkmate.
In the leftmost five files we have a pawn game
with initial position~[5], arranged symmetrically about the middle rank.
An ``immediate win'' in this game is a pawn that can promote to Queen
in three moves, ending the game by checkmate.
We may thus assume that, as in our pawn game of the previous section,
both sides play to prevent an ``immediate win''$\!$,
and the leftmost five files will eventually be empty or blocked.
This is why we have chosen a chessboard with an odd number of rows:
with an even number, as on the orthodox $8\times 8$ board,
one side's pawns would be at least one move closer to promotion,
and we would have to work harder to find positions in which,
as in Diags.~1A and~1B, an ``immediate win'' in the pawns game
by either player yields a chess win for the same player.
Once play ends in the [5]~component, we see why the component
in the g- and h-file is needed: the Zugzwangs arising from the
$[m]$ components all end with blocked pawns,
and if those were the only components on the board
then the chess game would end in stalemate,
regardless of which side ``won'' the pawns game.
But, with the g- and h-files on the board,
the side who lost the pawns game
must move in the central Zugzwang and lose the chess game.

To see how this happens, suppose that White is to move in Diag.~6A.
White must start 1~c5, the only winning move by the analysis in the
previous section.  Play may continue dxc5 2~dxc5 b5
(Black is lost, so tries to confuse matters with a loony move)
3~axb5 (declining the rope 3~a5?\ e5) axb5 4~e5 (Diag.~6B).
With all other pawns blocked, Black must now play 4\ldots g6
5~hxg6 g5.  If now 6~h7?\ g4 ends in stalemate, so White first plays
7~gxh5 (or even 7~h4), and then promotes the pawn on g6 and gives
checkmate in three more moves.

This construction clearly generalizes to show that any instance
of our pawn game supported on a board of length~$n$ can be realized
by a King-and-pawn endgame on any chessboard of at least $n+6$ files
whose height is an odd number greater than~$5$.

\vspace*{3ex}

\centerline{ {\large\bf 5.~Stopped files} }

By embedding our pawn game into generalized chess, we have constructed
a new class of endgames that can be analyzed by combinatorial game
theory.  But we have still not attained our aim of finding higher
Nimbers, because by Thm.~2 all the components of our endgames have
value~$0$ or~$*1$.  To reach $*2$ and beyond, we modify our components
by {\em stopping} some files.  We illustrate with Diag.~7A:

\vspace*{2ex}

\centerline{
\def\caption{Diag.~7A: whoever moves wins}
\longboardplus
{...........k}
{p.........pP}
{P.....p...P.}
{pppp...p.p..}
{......pP....}
{PPPP..P..P..}
{p......P..p.}
{P.........Pp}
{...........K}
\quad
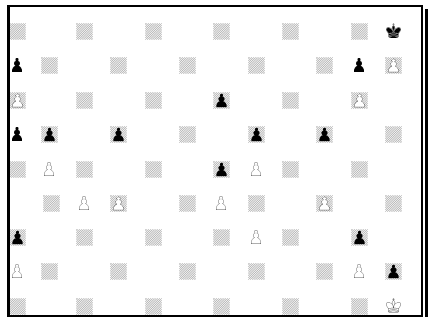
\def\caption{Diag.~7B: after 1 b5 cxb5 2 axb5}
\longboardplus
{...........k}
{p.........pP}
{P.....p...P.}
{pp.p...p.p..}
{.P....pP....}
{..PP..P..P..}
{p......P..p.}
{P.........Pp}
{...........K}
\
}

We have replaced the component~[5] of Diag.~6A by two components.
One is familiar: on the j-file (third from the right) we see
$[1]\cong *1$.  On the leftmost four files we have a new configuration.
This component looks like~[4], but with four extra pawns on the a-file.
These pawns are immobile, but have the effect of stopping the file
on both sides, so that a White pawn reaching a6 or a Black pawn
reaching a4 can no longer promote.  Without these extra pawns,
Diag.~7A would evaluate to $[4]+[1]\cong *1+*1=0$ and would thus
be a mutual Zugzwang.  But Diag.~7A is a first-player win, with the
unique winning move b5!. Indeed, suppose White plays b5 from Diag.~7A.
If Black responds 1\ldots axb5 then White's reply 2~axb5 produces
[:2]+[1] and wins.  So Black instead plays 1\ldots cxb5, expecting
the loony reply 2~cxb5.  But thanks to the stopped d-file
White can improve with 2~axb5!.  See Diag.~7B.  If now 2\ldots a5,
this pawn will get no further than~a4, while White forces a winning
breakthrough with 3~c5, for instance 3\ldots dxc5 4~dxc5 bxc5 5~b6 c4
6~b7 (Diag.~7C) and mates in two.  Notice that the extra pawns on
the a-file do not stop the b-file: if Black now captures the pawn
on~b7 then the pawn on~a7 will march in its stead.  We conclude that
in Diag.~7B Black has nothing better than 2\ldots axb5 3~cxb5, which
yields the same lost position \hbox{($[1]+[1]\cong0$)}
that would result from 2\ldots axb5.
After 1~b5 Black could also try the tricky 1\ldots c5, attempting
to exploit the a-file stoppage by sacrificing the \hbox{a6-pawn.}
After the forced 2~dxc5 (d5?\ axb5 3~axb5 j5 wins) dxc5 (Diag.~7D),
White would indeed lose after 3~bxa6?\ j5, but the pretty 3~a5!~wins.
However Black replies, a White pawn will next advance or capture to~b6,
and three moves later White will promote first and checkmate~Black.

%\vspace*{2ex}

\pagebreak

\centerline{
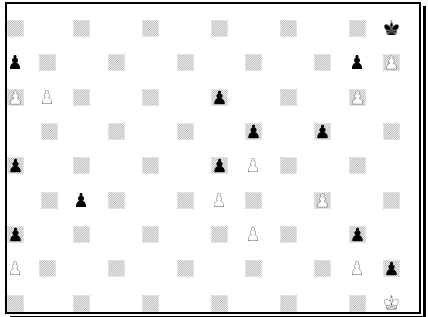
\def\caption{Diag.~7C: after 6 b7}
\longboardplus
{...........k}
{p.........pP}
{PP....p...P.}
{.......p.p..}
{p.....pP....}
{..p...P..P..}
{p......P..p.}
{P.........Pp}
{...........K}
\quad
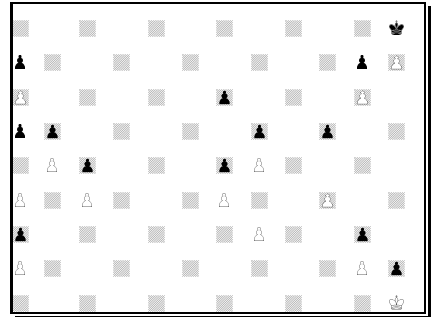
\def\caption{Diag.~7D: after 1\ldots c5 2~dxc5 dxc5}
\longboardplus
{...........k}
{p.........pP}
{P.....p...P.}
{pp.....p.p..}
{.Pp...pP....}
{P.P...P..P..}
{p......P..p.}
{P.........Pp}
{...........K}
\
}

Diag.~7A remains a first-player win even without the [1]~component
on the j-file (Diag.~7E).
The first move d5 wins as in our analysis of~[4]:
either cxd5 cxd5 or c5 bxc5 bxc5 a5 produces a decisive Zugzwang.
In fact, d5 is the only winning move in Diag.~7E.
The move c5 is loony as before (bxc5 bxc5 d5/dxc5).
With the a-file stopped, a5 is loony as well.
The opponent will answer b5 (Diag.~7F),
and if then cxb5 axb5!, followed by c5 and wins
while the pawn left on a5 is useless as in Diag.~7B.
This leaves (from Diag.~7F) c5, again producing the loony [:1].
Thus a5 is itself loony as claimed.

\vspace*{2ex}

\centerline{
\def\caption{Diag.~7E: whoever moves wins}
\longboardplus
{...........k}
{p.........pP}
{P.....p...P.}
{pppp...p....}
{......pP....}
{PPPP..P.....}
{p......P..p.}
{P.........Pp}
{...........K}
\quad
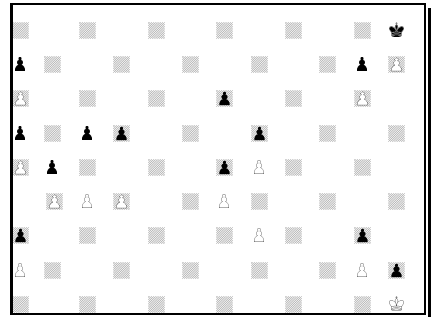
\def\caption{Diag.~7F: after 1 a5 b5}
\longboardplus
{...........k}
{p.........pP}
{P.....p...P.}
{p.pp...p....}
{Pp....pP....}
{.PPP..P.....}
{p......P..p.}
{P.........Pp}
{...........K}
\
}

Therefore the component in files a--d of Diag.~7A and Diag.~7E
is equivalent to an impartial game in which either side
may move to either~$0$ (with d5) or~$*1$ (with b5).
Hence this component has the value $*2$!
On longer boards of odd height $\geq9$,
we can stop some of the files in $[m]$ for other~$m$.
We next show that each of the resulting components
is equivalent to a Nim-heap, some with values $*4$, $*8$ and beyond.

%\vspace*{3ex}

\pagebreak

\centerline{ {\large\bf 6.~The pawns game with stopped files} }

{\bf 6.1: Game definition and components}

We modify our pawn game by choosing a subset of the $n$ files
and declaring that the files in that subset are ``stopped''$\!$.
A pawn reaching its opposite row now scores an immediate win
only if it is on an unstopped file.\footnote{
  If the file is stopped, the pawn does not ``promote'':
  it remains a pawn, and can make no further moves.
  Recall that this was the fate of Black's a5-pawn in Diag.~7C.
  }
{\em We require that no two stopped files be adjacent.}
This requirement arises naturally from our implementation
of stopped files in King-and-pawn endgames on large chessboards.
As it happens, the requirement is also necessary
for our analysis of the modified pawns game.  For instance, if 
adjacent stopped files were allowed then a threat to capture a pawn
might not be an entailing move.

In the last section we already saw the effects of stopped files on 
the play of the game.  We next codify our observations.  We begin
by extending our notation for quiescent components.  In \S3, such
a component was entirely described by the number~$m$ of consecutive
initial files that the component comprises.  In the modified game,
we must also indicate which if any of these $m$ files is stopped.
We denote a stopped initial file by~$1$, and an unstopped one by~$0$.
A string of $m$ binary digits, with no $1$'s adjacent, then denotes
a quiescent component of $m$ initial files.  For instance, the
component we called~$[m]$, with no stopped files, now becomes
$[000\cdots0]=[0^m]$; the component with value $*2$ on files a--d
of Diag.~7A is denoted $[1000]$.
An initial file that may be either open or closed
will be denoted $i$ (or $i\0_1$, $i\0_2$, etc.);
an arbitrary ``word'' of $0$'s and $1$'s
will be denoted $w$ (or $w'$, $w\0_1$, $w\0_2$, etc.).

A component comprising just one initial file, stopped or not, still
has value $\{0|0\}=*1$.  In a component of at least two initial files,
every move threatens to capture and is entailing.  This is true even
if the file(s) of the threatened pawn(s) is or are stopped, because
an immediate win is then still threatened by advancing in that stopped
file, as we saw in Diag.~7D where White wins by 3~a5!.  (Here
we need the condition that two adjacent files cannot both be stopped:
if in Diag.~7D the b-file were somehow stopped as well then 3~a5 would
lose to either 3\ldots bxa5 or 3\ldots axb5.)

Consider first a move by the pawn on the first or last file
of the component (without loss of generality: the first), attacking
just one pawn.  As in \S3, the opponent must move the attacked pawn
on the second file, either advancing it or capturing the attacking pawn.
In the latter case, the pawn must be re-captured, and the sequence has
the effect of removing the component's first two files.
In the latter case, the component becomes quiescent if it had only two
files (in which case the first move in the component was loony,
as before); otherwise the advanced pawn in turn attacks a third-file
pawn, which must capture or advance.  But now a new consideration
enters: if the {\em first}\/ file is stopped, then the capture loses
immediately since the opponent will re-capture from the first file
and touch down on the second, necessarily unstopped, file.
(See Diag.~7F after 2~cxb5 axb5.)  Note that the first file,
though closed, can still affect play for one turn after its closure
if it is stopped.  We thus need a notation for such files,
as well as stopped ``:''~files, which may become closed.
We use an underline: a stopped ``:''~file will be denoted ``$\ul:$'',
and a stopped blocked file will simply be denoted~``$\sb$''.
Thus the moves discussed in this paragraph are as follows,
with each $w$ denoting an arbitrary word {\em of positive length}\/:
\begin{itemize}
\item From [0] or [1], either side may move to~$0$.
\item From $[0w]$ or $[1w]$, either side may move to
 [:$w$] or [\ul{:}$w$] respectively.
\item a move to [:0], [:1], or [\ul{:}0] is loony.
\item {}[:0$w$] entails a move to [:.]+[$w$] or [:$w$];
 [\ul{:}0$w$] entails a move to [\ul{:}.]+[$w$] or [\sb:$w$];
 and [:1$w$] entails a move to [:.]+[$w$] or [\ul{:}$w$].
 Each of [:.] and [\ul{:}.] entails a move to~$0$.
\item {}[\sb:0] or [\sb:1] entails a move to~$0$;
 and [\sb:0$w$] or [\sb:1$w$] entails a move to
 [:$w$] or [\ul{:}$w$] respectively.
\end{itemize}

It remains to consider a pawn move in the interior of a quiescent
component.  Such a move attacks two of the opponent's pawns, and
entails a capture.  If neither of the attacked pawns is on a stopped
file, then either of them may capture, forcing a re-capture from
the same file, just as in the pawn game without stopped files.
If both pawns are on stopped files (see Diag.~8A), then a capture from
either of these files can be met by a capture from the other file
(Diag.~8B), forcing a further capture and re-capture to avoid
immediate loss.  The first player may also choose to make
the first re-capture from the same file (Diag.~8C),
but we can ignore this possibility because the opponent can
still re-capture again to produce the same position as before,
but has the extra option of advancing the attacked pawn.

\vspace*{2ex}

\centerline
{
\boardthree
{pppppppp}
{...P....}
{PPP.PPPP}
{\ Diag.~8A (c, e files blocked)}
\
\boardthree
{pppp.ppp}
{...P....}
{PP..PPPP}
{Diag.~8B} % after 1\ldots exd2 2 cxd2
\
\boardthree
{pppp.ppp}
{...P....}
{PPP..PPP}
{Diag.~8C} % after 2 exd2(?) (2\ldots cxd2, also c2)
\
}

Finally,
if just one of the two attacked pawns lies on an unstopped file,
it may as well capture, forcing a re-capture in the same file:
capturing with the other pawn lets the first player capture
from the stopped file, forcing a further capture and re-capture
to avoid immediate loss, and thus denying the opponent the option
to capture with one attacked pawn and then advance the other.
We next summarize the moves discussed in this paragraph
that we did not list before.  Here $w,w\0_1,w\0_2$ again denote
arbitrary words, which may be empty (length zero) 
except for the first item:

\begin{itemize}
\item From $[w\0_10w\0_2]$ or $[w\0_11w\0_2]$
 with $w\0_1,w\0_2$ of positive length, either side may move to
 [$w\0_1$:$w\0_2$] or [$w\0_1$\ul{:}$w\0_2$] respectively.
\item {}[$w\0_1 i\0_1$:$i\0_2 w\0_2$] entails a move
 to [$w\0_1$]+[.:$i\0_2 w\0_2$] or [$w\0_1 i\0_1$:.]+[$w\0_2$];
 likewise, [$w\0_1010w\0_2$] entails a move to
 [$w\0_10$]+[.\ul{:}$w\0_2$] or [$w\0_1$\ul{:}.]+[$0w\0_2$].
\item {}[.:$w$] entails a move to [:$w$];
 likewise, [.\ul{:}$w$] entails a move to [\ul{:}$w$].
\item A move to [$w\0_1$1:1$w\0_2$] is equivalent to a move to
 $[w\0_1]+[w\0_2]$.
\item A move to [$w\0_1$1:0$w\0_2$] is equivalent to a move to
 $[w\0_1]+[:\!0w\0_2]$.
\end{itemize}
Only the last two cases are directly affected by stopped files.

\vspace*{2ex}

{\bf 6.2: Reduction to Nim}

Even with stopped files it turns out that our pawn game 
still reduces to an impartial game, and thus to Nim, once
immediately losing and loony moves are eliminated:

{\bf Theorem 3.}
{\sl
i) Each component $[w]$ is equivalent to a Nim-heap
 of some size~$\epsilon(w)$.\\
ii) A move to $[:\!i]$, or $[\ul{:}0i]$ is loony.
 For each $w$ of positive length, a move to $[:\!iw]$
 is either loony or equivalent to a move to $[w]\cong *\epsilon(w)$.
 The move to $[:\!0w]$ or $[:\!1w]$ is loony
 if and only if a move to $[:\!w]$ or $[\ul{:}w]$ respectively
 is \ul{not} loony and is equivalent to a move to $*\epsilon(w)$.
 For each $w$ of positive length, a move to $[\ul{:}0iw]$
 is either loony or equivalent to a move to $[iw]\cong *\epsilon(iw)$.
 The move to $[\ul{:}00w]$ or $[\ul{:}01w]$ is \ul{not} loony
 if and only if a move to $[:\!w]$ or $[\ul{:}w]$ respectively
 is \ul{not} loony and is equivalent to a move to $*(0w)$ or $*(1w)$.\\
iii) For any words $w\0_1,w\0_2$, a move to $[w\0_10\!:\!0w\0_2]$
 or $[w\0_10\ul{:}0w\0_2]$ is either loony or equivalent to a move to
 $[w\0_1]+[w\0_2] = *(\epsilon(w\0_1) \nimsum \epsilon(w\0_2))$.
 The move to $[w\0_10\!:\!0w\0_2]$ is non-loony if and only both~$w\0_1$
 and $w\0_2$ satisfy the criterion of~(ii) for a move to~$[:\!0w]$
 to be non-loony.  Likewise, the move to $[w\0_10\ul{:}0w\0_2]$
 is non-loony if and only both~$w\0_1$ and $w\0_2$ satisfy the criterion
 of~(ii) for a move to~$[\ul{:}0w]$ to be non-loony.\\
iv) The function $\epsilon$ from strings of $0$'s and $1$'s with
 no consecutive $1$'s to nonnegative integers is recursively determined
 by (ii) and (iii): $\epsilon(i)=1$, and for $w$ of length~$>1$
 the value $\epsilon(w)$ is the mex of the values of the Nim
 equivalents of all non-loony moves as described in (ii),~(iii).
}

This is proved in exactly the same way as Thm.~1.  Note that we
do not evaluate moves to [\sb:$w$].  Such a move is available only
if the opponent just moved to [\ul{:}0$w$].  If that move was loony
then we win, capturing unless the other components sum to $\epsilon(w)$
in which case we advance, forcing the opponent to advance in return
and winning whether that forced advance was loony or not.  If the
opponent's move to [\ul{:}0$w$] was not loony then capturing or
advancing our attacked pawn yields equivalent positions.~~$@p$

{\bf 6.3 Numerical results}

Thm.~3 yields a practical algorithm for evaluating $\epsilon(w)$.
If $w$ has length $m$, the recursion in~(iv) requires $O(m^2 \log m)$
space, to store $\epsilon(w')$ as it is computed for each substring~$w'$
of consecutive bits of~$w$, and $O(m^3)$ table lookups and nim-sums,
to recall each $\epsilon(w')$ as it is needed and combine pairs.
Unlike the situation for the simple game with no files stopped,
where we obtained a simple closed form for $\epsilon(m)$ (Thm.~2),
here we do not know such a closed form.  We can, however implement
the $O(m^3)$ algorithm to compute $\epsilon(w)$ for many~$w$.
We conclude this paper with a report on the results of several
such computations and our reasons for believing that $\epsilon(w)$
can be arbitrarily large.

We saw already that $\epsilon(1000)=2$; this is the unique $w$ of
minimal length such that $[w]$ has value $*2$, except that
the reversal $[0001]$ of $[1000]$ has the same value.  Clearly
$[0]$ and $[1]$ are the smallest instances of~$*1$.
We first find $*4$, $*8$ and $*16$ at lengths 9, 20, and 43,
for $w=101001000$, 10100100010100001000,
and 1010010001000000010100010000000101000100101, among others.
The following table lists for each $k\leq16$ the least $m$ such that
$\epsilon(w)=k$ for some word of length~$m$:

\centerline{
\begin{tabular}{c|cccccccccccccccc}
$k$ & 1 & 2 & 3 & 4 & 5 & 6 & 7 & 8 & 9 &10 &11 &12 &13 &14 &15 &16
\\ \hline
$m$ & 1 & 4 & 6 & 9 &11 &14 &16 &20 &22 &25 &27 &30 &32 &37 &39 &43
\end{tabular}
}

It seems that, for each $k$, instances of $*k$\/ are quite
plentiful as $m$ grows.  The following table gives for each
$35\leq m \leq 42$ the proportion of length-$m$ words~$w$
with $w(m)=0,1,2,\ldots,9$, rounded to two significant figures:\footnote{
  To compute such a table one need do $O(m)$ basic operations for each
  choice of~$w$, rather than $O(m^3)$, because $\epsilon(w')$ has already
  been computed for each substring~$w'$ --- this as long as one has enough
  space to store $\epsilon(w')$ for all $w'$ of length at most $m-2$.
  }

\centerline{
\begin{tabular}{c|cccccccccc}
$m$ &  0   &  1   &  2   &  3   &   4   &   5   &   6   &   7   &  8   &  9
\\ \hline
35  & 24\% & 26\% & 19\% & 15\% & 5.4\% & 5.7\% & 2.7\% & 2.5\% &.51\% &.25\%
\\
36  & 22\% & 27\% & 18\% & 15\% & 5.5\% & 5.7\% & 2.6\% & 2.8\% &.54\% &.27\%
\\
37  & 26\% & 22\% & 14\% & 19\% & 5.8\% & 5.5\% & 2.8\% & 2.8\% &.55\% &.31\%
\\
38  & 25\% & 23\% & 16\% & 17\% & 5.7\% & 5.7\% & 3.1\% & 2.7\% &.56\% &.35\%
\\
39  & 22\% & 26\% & 19\% & 14\% & 5.6\% & 5.9\% & 3.0\% & 3.0\% &.59\% &.37\%
\\
40  & 24\% & 24\% & 16\% & 18\% & 5.9\% & 5.7\% & 3.0\% & 3.2\% &.61\% &.40\%
\\
41  & 26\% & 22\% & 15\% & 19\% & 5.9\% & 5.8\% & 3.3\% & 3.1\% &.61\% &.44\%
\\
42  & 22\% & 24\% & 18\% & 15\% & 5.8\% & 6.0\% & 3.3\% & 3.2\% &.63\% &.47\%
\end{tabular}
}

Especially for $*0$ through $*3$, the proportions seem to be bounded
away from zero but varying quite erratically with~$m$.  The small
proportions of $*6$ through~$*9$ appear to rise slowly but not smoothly.
We are led to guess that for each $k$\/ there are length-$m$ components
of value $*k$\/ once $m$ is large enough --- perhaps $m\gg k$\/ suffices
--- and ask for a description and explanation of the proportion of
components of value~$*k$\/ among all components of length~$m$.
In particular, is it true for each~$k$\/ that this proportion
is bounded away from zero as $m\mra\infty$?

It is well known that the number of binary words of length~$m$ without
two consecutive $1$'s is the $(m+2)$-nd Fibonacci number.  This number
grows exponentially with~$m$, soon putting an end to exhaustive
computation.  We do not expect to be able to extend such computations
to find the first $*32$, which probably occurs around $m=90$.
Nevertheless we have reached $*32$ and much more
by computing $\epsilon(w)$ for periodic $w$ of small period~$p$.
This has the computational advantage that for each $m'<m$ we have
at most $p$ substrings of length~$m'$ to evaluate, rather than
the usual $m+1-m'$.\footnote{
  Once this is done for some period-$p$ pattern, one can also
  efficiently evaluate components with the same repeated pattern
  attached to any initial configuration of blocked and unblocked files.
  We have not yet systematically implemented this generalization.
  }

We have done this for various small~$p$.  Often the resulting
Nim-values settle into a repeating pattern, of period~$p$ or
some multiple of~$p$.  This is what happened in Thm.~2 for $p=1$,
with period~$10p$.  Usually the multiplier is smaller than~$10$,
though blocking files $14r$ and $14r+5$ produces a period of
$504=36\cdot 14$.

All repeating patterns with $p\geq 4$ soon become periodic,
but for larger~$p$ some choices of pattern yield large and
apparently chaotic Nim-values.  For instance, we have reached
$*4096$ by blocking every sixth file in components of length
up to $2\cdot 10^5$.  For each $\alpha=3,4,\ldots,12$ the shortest
such component of value $*(2^\alpha)$ has files $6r+4$ blocked,
with length~$n$ given by the following table:

\centerline{
\begin{tabular}{c|cccccccccc}
$\alpha$ &  3  &  4  &  5  &  6  &   7  &   8  &   9  &   10  &   11  &   12
\\ \hline
  $n$    & 51  & 111 & 202 & 497 & 1414 & 3545 & 8255 & 21208 & 61985 & 187193
\end{tabular}
}

This again suggests that all $*k$\/ arise: even if the Nim-values
for $p=6$ ultimately become periodic, we can probably re-introduce
chaos by blocking or unblocking a few files.  Of course we have no
idea how to prove that arbitrarily large~$k$\/ appear this way.

Finally, for a few repeating patterns we observe behavior apparently
intermediate between periodicity and total chaos.  Blocking every
fifth or tenth file yields Nim-values that show some regularity
without (yet?)\ settling into a period.  Indeed in both cases
the largest values grow as far as we have extended the search
(through length $10^5$), though more slowly and perhaps
more smoothly than for $p=6$.  Such families of components
seem the most likely place to find and prove an arithmetic
periodicity or some more complicated pattern that finally proves
that all $*k$\/ arise and thus fully embeds Nim into generalized
King-and-pawn endgames.

\vspace*{3ex}

\centerline{ {\large\bf Acknowledgements} }
This paper was typeset in \LaTeX, using Piet Tutelaers' chess
fonts for the Diagrams.  The research was made possible in part
by funding from the Packard Foundation.  I thank the Mathematical
Sciences Research Institute for its hospitality during the completion
of this paper.

\end{document}

%% file: macros2.tex
\newlength{\squarewidth} 
\newlength{\diagwidth} 
\font\shmt=chess13+
\font\cfig=chess10
\catcode`@=\active
\def@#1{\lower1.3pt\hbox{\cfig\if#1SN\else#1\fi}}
\def\sqskip{\hspace*\squarewidth}
\def\ltod#1{
  \if#1KJ\else
  \if#1kj\else
  \if#1QL\else
  \if#1ql\else
  \if#1RS\else
  \if#1rs\else
  \if#1BA\else
  \if#1ba\else
  \if#1NM\else
  \if#1nm\else
  \if#1SM\else
  \if#1sm\else
  \if#1PO\else
  \if#1po\else
  Z
  \fi\fi\fi\fi\fi\fi\fi\fi\fi\fi\fi\fi\fi\fi
  }
\def\llhack#1{                %to make {\shmt 0} (empty light square)
  \if#10\sqskip\else          %act the way it should, and make S=N
  \if#1.\sqskip\else
  \if#1SN\else\if#1sn\else#1\fi\fi\fi\fi}
\def\evenrank#1#2#3#4#5#6#7#8{
  \llhack#1\ltod#2\llhack#3\ltod#4\llhack#5\ltod#6\llhack#7\ltod#8}
\def\oddrank#1#2#3#4#5#6#7#8{
  \ltod#1\llhack#2\ltod#3\llhack#4\ltod#5\llhack#6\ltod#7\llhack#8}
\def\board#1#2#3#4#5#6#7#8#9{
  \setlength{\unitlength}{1\squarewidth}
  \begin{picture}(9.5,9.8)
  {\shmt
  \put(.75,8.52){\evenrank #1}
  \put(.75,7.52){\oddrank  #2}
  \put(.75,6.52){\evenrank #3}
  \put(.75,5.52){\oddrank  #4}
  \put(.75,4.52){\evenrank #5}
  \put(.75,3.52){\oddrank  #6}
  \put(.75,2.52){\evenrank #7}
  \put(.75,1.52){\oddrank  #8}}
  \thinlines
  \put(.75,1.5){\framebox(8,8){}}
  \thicklines
  \put(.85,1.35){\line(1,0){8.07}}\put(8.9,9.4){\line(0,-1){8.07}}
  \put(.75,0.5){\makebox(8,1)[b]{#9}}
  \end{picture}
  }
\def\longevenrank#1#2#3#4#5#6#7#8#9{
  \llhack#1\ltod#2\llhack#3\ltod#4\llhack#5\ltod#6\llhack#7\ltod#8\evenrest#9}
\def\evenrest#1#2#3#4{\llhack#1\ltod#2\llhack#3\ltod#4}
\def\longoddrank#1#2#3#4#5#6#7#8#9{
  \ltod#1\llhack#2\ltod#3\llhack#4\ltod#5\llhack#6\ltod#7\llhack#8\oddrest#9}
\def\oddrest#1#2#3#4{\ltod#1\llhack#2\ltod#3\llhack#4}

\def\longboardplus#1#2#3#4#5#6#7#8#9{
  \setlength{\unitlength}{1\squarewidth}
  \begin{picture}(13.5,10.8)
  {\shmt
  \put(.75,9.52){\longoddrank #1}
  \put(.75,8.52){\longevenrank #2}
  \put(.75,7.52){\longoddrank  #3}
  \put(.75,6.52){\longevenrank #4}
  \put(.75,5.52){\longoddrank  #5}
  \put(.75,4.52){\longevenrank #6}
  \put(.75,3.52){\longoddrank  #7}
  \put(.75,2.52){\longevenrank #8}
  \put(.75,1.52){\longoddrank  #9}}
  \thinlines
  \put(.75,1.5){\framebox(12,9){}}
  \thicklines
  \put(.85,1.35){\line(1,0){12.07}}\put(12.9,10.4){\line(0,-1){9.07}}
  \put(.75,0.5){\makebox(12,1)[b]{\caption}}
  \end{picture}
  }
\def\boardthree#1#2#3#4{
  \setlength{\unitlength}{1\squarewidth}
  \begin{picture}(9.5,4.8)
  {\shmt
  \put(.75,3.52){\oddrank  #1}
  \put(.75,2.52){\evenrank #2}
  \put(.75,1.52){\oddrank  #3}}
  \thinlines
  \put(.75,1.5){\framebox(8,3){}}
  \thicklines
  \put(.85,1.35){\line(1,0){8.07}}\put(8.9,4.4){\line(0,-1){3.07}}
  \put(.75,0.5){\makebox(8,1)[b]{#4}}
  \end{picture}
  }
\def\longboardthree#1#2#3#4{
  \setlength{\unitlength}{1\squarewidth}
  \begin{picture}(13.5,4.8)
  {\shmt
  \put(.75,3.52){\longoddrank  #1}
  \put(.75,2.52){\longevenrank #2}
  \put(.75,1.52){\longoddrank  #3}}
  \thinlines
  \put(.75,1.5){\framebox(12,3){}}
  \thicklines
  \put(.85,1.35){\line(1,0){12.07}}\put(12.9,4.4){\line(0,-1){3.07}}
  \put(.75,0.5){\makebox(12,1)[b]{#4}}
  \end{picture}
  }